\input amstex 
\documentstyle{amsppt} 
\loadbold
\magnification=1200
\NoBlackBoxes
\pagewidth{6.4truein}
\pageheight{8.5truein}

\define\C{{\Bbb C}}
\redefine\D{{\Bbb D}}
\define\N{{\Bbb N}}
\define\R{{\Bbb R}}
\define\T{{\Bbb T}}
\define\A{{\Cal A}}

\define\F{{\Cal F}}
\redefine\O{{\Cal O}}
\redefine\d{{\partial}}
\redefine\phi{{\varphi}}
\redefine\epsilon{{\varepsilon}}
\define\ac{\acuteaccent}
\define\gr{\graveaccent}
\define\set#1{\{#1\}}
\define\PSH{\operatorname{PSH}}
\define\pr{\operatorname{pr}}

\redefine\cdot{\boldsymbol\cdot}

\hyphenation{pluri-sub-har-mon-ic}

\refstyle{B}
\NoRunningHeads
\TagsOnRight

\topmatter 
\title Plurisubharmonic extremal functions, \\
Lelong numbers and coherent ideal sheaves
\endtitle
\author Finnur L\ac arusson and Ragnar Sigurdsson \endauthor
\address Department of Mathematics, University of Western Ontario,
	London, Ontario N6A~5B7, Canada \endaddress
\email larusson\@uwo.ca \endemail
\address  Department of Physics and Mathematics, 
Mid Sweden University, S-851 70 Sundsvall, Sweden \endaddress
\email ragnar\@fmi.mh.se \endemail 
\date 14 July 1998  \enddate 

\thanks The first-named author was supported in part by the Natural
Sciences and Engineering Research Council of Canada.
\endthanks

\abstract We introduce a new type of pluricomplex Green function which
has a logarithmic pole along a complex subspace $A$ of a complex
manifold $X$.  It is the largest negative plurisubharmonic function on
$X$ whose Lelong number is at least the Lelong number of
$\log\max\set{|f_1|,\dots,|f_m|}$, where $f_1,\dots,f_m$ are local
generators for the ideal sheaf of $A$.  The pluricomplex Green function
with a single logarithmic pole or a finite number of weighted poles is a
very special case of our construction.  We give several equivalent
definitions of this function and study its properties, including
boundary behaviour, continuity, and uniqueness.  This is based on and
extends our previous work on disc functionals and their envelopes. 
\endabstract



\endtopmatter

\document

\specialhead 1. Introduction
\endspecialhead

\noindent
Let $X$ be a complex manifold.  For each function ${\alpha}:X\to
[0,+{\infty})$ we let
 $$
\F_\alpha=\set{u\in \PSH(X)\, ;\, u\leq 0, {\nu}_u\geq {\alpha}}
 $$
and
 $$
G_{\alpha}=\sup\F_\alpha,
 $$
where $\PSH(X)$ is the class of plurisubharmonic functions on $X$
(including the constant function $-\infty$) and ${\nu}_u$ denotes the
Lelong number of $u$.  Then
$G_\alpha$ is plurisubharmonic and 
$G_\alpha\in\F_\alpha$ \cite{6, Prop\. 5.1}.  Recall that the Lelong number
is a biholomorphic invariant, and if $u$ is plurisubharmonic in a 
neighbourhood of $0$ in $\C^n$, then
$$
\nu_u(0)=\lim_{r\to 0}\frac{\sup\limits_{|z|=r} u(z)}{\log r}.
$$

If ${\alpha}$ is the characteristic function of a one-point set $\set
a$, then $G_{\alpha}$ is the pluricomplex Green function $G_a$ of $X$
with a logarithmic pole at $a$, first defined by Klimek \cite{4}.  Such
functions have also been studied e.g\.  by Demailly \cite{1}, Edigarian
\cite{2, 3}, Lempert \cite{9, 10}, and Zeriahi \cite{17}.  If ${\alpha}$
has finite support, then $G_{\alpha}$ is the pluricomplex Green function
of $X$ with a logarithmic pole of weight $\alpha(a)$ at each point $a$
of the support.  Such functions were first defined by Lelong \cite{8}. 

This paper is a study of the much larger class of functions
$G_A=G_\alpha$, where $A$ is a (closed) complex subspace of $X$, and
$\alpha=\nu_A$ is the Lelong number of the plurisubharmonic function
$\log\max\set{|f_1|,\dots,|f_m|}$, where $f_1,\dots,f_m$ are local
generators for the ideal sheaf of $A$.  This number is independent of
the choice of generators.  We call $G_A$ the {\it pluricomplex Green
function with a logarithmic pole along} $A$, or simply the {\it Green
function with a pole along} $A$.

Our results here are based on our previous work \cite{6} in the theory
of disc functionals and their envelopes, which in turn builds on the
pioneering work of Poletsky \cite{11, 12}.  Indeed, on domains in Stein
manifolds (and in fact under certain much weaker conditions), $G_A$ is
the envelope of the Lelong functional associated to the function
$\nu_A$.  In Section 2 we review the theory of disc functionals with
emphasis on the three known classes of examples.  We add a few new
results to the theory of the Lelong functional, including a product
property, generalizing a result of Edigarian \cite{3}. 

In Section 3 we study the pluricomplex Green function with a logarithmic
pole along a complex subspace.  Our main results may be summarized as
follows. 

\proclaim{Main Theorem} Let $X$ be a relatively compact domain in a
Stein manifold $Y$, and let $A$ be the intersection with $X$ of a
complex subspace $B$ of $Y$.  Then $G_A$ is locally bounded and maximal
on $X\setminus A$, and $\nu_{G_A}=\nu_A$. 

If $X$ has a strong plurisubharmonic barrier at $p\in\d X\setminus B$,
then $G_A(x)\to 0$ as $x\to p$. 

If $A$ is a divisor, then the Levi form of $G_A$ is at least $\pi$ times
the current of integration over $A$. 

If $X$ has a strong plurisubharmonic barrier at every boundary point and
$B$ is a principal divisor, then the set of points in $X$ at which $G_A$
is discontinuous is pluripolar.
\endproclaim

We relate $G_A$ to the Poisson and Riesz functionals.  This yields several
alternative definitions of $G_A$.  We show that $G_A$ is uniquely
determined by some of its key properties when $A$ is a principal
divisor.  We also present a few instructive examples.  Among other
things, through an investigation of Green functions, we obtain a bounded
pseudoconvex domain in $\C^2$ such that $\log\tanh$ of the Carath\ac
eodory distance to the boundary is not plurisubharmonic. 

\smallskip\noindent
{\it Acknowledgements.}  We would like to thank Zbigniew B\l ocki,
Evgeny Poletsky, Azim Sadullaev, and Ahmed Zeriahi for valuable
discussions. 

\vfill\eject

\specialhead  2. Disc functionals and their envelopes
\endspecialhead

\noindent
Let $X$ be a complex manifold. A holomorphic map from the unit disc
$\D$ to $X$ is called an {\it analytic disc} in $X$, and if it can be
extended holomorphically to some neighbourhood of
$\overline \D$ then it is called a {\it closed analytic disc} in $X$.
We let $\O(\D,X)$ denote the set of all analytic discs in $X$ and
$\A_X$ denote the set of all closed analytic discs in $X$.
A {\it disc functional} on $X$ is a map
$H$ from a subset of $\O(\D,X)$ containing
$\A_X$ to $\R\cup \set{\pm{\infty}}$. 
The {\it envelope} of $H$ is the
function $EH:X\to \R\cup \set{\pm{\infty}}$ defined by the formula
 $$
EH(x)=\inf\set{H(f)\, ;\, f\in \A_X, f(0)=x}.
 $$ 
Envelopes of disc functionals were first defined and studied by
Poletsky \cite{11, 12}.

In our paper \cite{6}, we studied three classes of examples of disc
functionals and proved that their envelopes are plurisubharmonic for
a large collection of manifolds $X$.  These functionals are called
the Poisson, Riesz, and Lelong functionals.

If ${\varphi}:X\to \R\cup \set{-{\infty}}$ is upper semicontinuous,
then the {\it Poisson functional} $H_P^{\varphi}$ is defined on
$\A_X$ by
$$H_P^{\varphi}(f)=\dfrac 1{2\pi} \int_\T\phi\circ f\,d\lambda,
$$
where $\lambda$ is the arc length measure on the unit circle $\T$.   

We define 
$$\Cal F_P^{\varphi}=\set{u\in\PSH(X)\,;\,u\leq\phi}.
$$ 
Then $\sup \Cal F_P^{\varphi}$ is plurisubharmonic,
$\sup \Cal F_P^{\varphi}\leq EH_P^\phi$,
and equality holds if and only if $EH_P^\phi$ is plurisubharmonic
\cite{6, Prop\. 2.1}.  

If  $v$ is a plurisubharmonic function on $X$, then the {\it Riesz
functional} $H_R^v$ is defined on $\O(\D,X)$ by
$$H_R^v(f)=\dfrac 1{2\pi} \int_\D\log|\cdot|\,\Delta(v\circ f)
$$
if $f\in \O(\D,X)$ and $v\circ f$ is not identically $-\infty$, 
where $\Delta(v\circ f)$ is considered as a positive Borel measure on
$\D$.  If $f\in\O(\D,X)$ and $v\circ f=-\infty$, then we set $H_R^v(f)=0$.

We define 
$$
\Cal F_R^v=\set{u\in\PSH(X)\,;\, u\leq 0,\,\Cal L(u)\geq\Cal L(v)}.
$$ 
Here, $\Cal L(u)$ denotes the Levi form $i\d\bar\d u$ of $u$.  We set
$\Cal L(-\infty)=0$.  If $v:X\to\R$ is continuous, 
then $\sup \Cal F_R^v$ is
plurisubharmonic, $\sup \Cal F_R^v\leq EH_R^v$, and equality holds if
and only if $EH_R^v$ is plurisubharmonic \cite{6, Thm\. 4.4}.

If  $\alpha$ is a nonnegative function on $X$,  then the {\it Lelong 
functional} $H_L^\alpha$ is defined on $\O(\D,X)$ by the formula
$$
H_L^{\alpha}(f)=\sum_{z\in\D} \alpha(f(z))\,m_z(f)\,\log|z|.
$$
The sum, which may be uncountable,
is defined as the infimum of its finite partial sums.
Here, $m_z(f)$ denotes the multiplicity of $f$ at $z$,
defined in the following way.  If $f$ is constant, let $m_z(f)=\infty$.
If $f$ is nonconstant, let $(U,\zeta)$ be a coordinate neighbourhood
on $X$ with $\zeta(f(z))=0$.  Then there exists an integer
$m$ such that $\zeta(f(w))=(w-z)^m g(w)$
where $g:V\to \C^n$  is a map defined in a neighbourhood $V$ of $z$ with
$g(z)\neq 0$.
The number $m$, which is independent of the choice of local coordinates,
is the multiplicity of $f$ at $z$.  

We define  
$$
\Cal F_L^{\alpha}=\set{u\in\PSH(X)\,;\,u\leq 0, {\nu}_u\geq{\alpha}}.
$$
Then $G_{\alpha}=\sup\Cal F_L^{\alpha}$ is plurisubharmonic.
If $u\in \Cal F_L^{\alpha}$ and $f\in \O(\D,X)$, then 
$u(f(0))\leq H_L^{\alpha}(f)$, so $G_{\alpha}\leq EH_L^{\alpha}$.
The function $EH_L^{\alpha}$ is plurisubharmonic if and only if
$EH_L^{\alpha}\in \Cal F_L^{\alpha}$, and then
$G_{\alpha}=EH_L^{\alpha}$ \cite{6, Prop\. 5.1}.  

Clearly,
$$
EH_L^\alpha=\inf EH_L^\beta,
$$
where the infimum is taken over all functions $\beta$ with finite support
such that $0\leq\beta\leq\alpha$.  When $EH_L^\beta$ is plurisubharmonic,
it is a pluricomplex Green function with finitely many weighted poles.

Let us now take a closer look at the Lelong
functional.   For $f\in \O(\D,X)$ we define $f^*{\alpha}:\D\to
[0,+{\infty})$ by the formula
 $$
f^*{\alpha}(z)={\alpha}(f(z))m_z(f), \qquad z\in \D.
 $$
Then 
 $$
{\mu}_f^{\alpha}=2{\pi}\sum_{z\in \D} f^*{\alpha}(z){\delta}_z
 $$
is a well defined positive Borel measure on $\D$, where ${\delta}_z$ is the
Dirac measure at $z$.  We define
 $$
v_f^{\alpha}(z)=\int_\D G(z,\cdot)\, d{\mu}_f^{\alpha}
=\sum_{w\in \D} f^*{\alpha}(w)\log\bigg|\dfrac{z-w}{1-\bar w z}\bigg|,  
\qquad z\in \D,
 $$
where $G$ denotes the Green function of the unit disc,
$$
G(z,w)=\dfrac 1{2{\pi}}\log\bigg|\dfrac{z-w}{1-\bar w z}\bigg|.
$$
Then $v_f^{\alpha}$ is subharmonic in $\D$.  We have $v_f^{\alpha}\neq
-{\infty}$ if and only if
 $$
\sum_{z\in \D} f^*{\alpha}(z)(1-|z|)=\dfrac 1{2{\pi}} 
\int_\D (1-|\cdot|)\, d{\mu}_f^{\alpha}<{\infty}.
 $$
If $v_f^{\alpha}\neq -{\infty}$, then ${\mu}_f^{\alpha}$ has 
finite mass on compact sets, the sum which defines
${\mu}_f^{\alpha}$ is convergent in the sense of distributions,
and ${\Delta}v_f^{\alpha}={\mu}_f^{\alpha}$.

The Lelong number ${\nu}_v(z)$ of a subharmonic function $v$ on a
domain in $\C$ at a point $z$ is ${\Delta}v(\set z)/2{\pi}$, i.e., the
Riesz mass of $v$ at the point $z$ divided by $2{\pi}$.  We therefore
have ${\nu}_{v_f^{\alpha}}=f^*{\alpha}$,
$H_L^{\alpha}(f)=v_f^{\alpha}(0)$, and from the Riesz representation
formula we see that $v_f^{\alpha}$ is the largest negative subharmonic
function on $\D$ satisfying ${\nu}_v\geq f^*{\alpha}$.

Assume now that $f\in \O(\D,X)$ is an extremal disc for the Lelong
functional, i.e., $G_{\alpha}(f(0))=H_L^{\alpha}(f)=v_f^{\alpha}(0)$,
and that $G_{\alpha}(f(0))>-{\infty}$.  Then the function
$v=G_{\alpha}\circ f$ is subharmonic in $\D$, $v\leq 0$, and
${\nu}_v(z)\geq {\nu}_{G_{\alpha}}(f(z))m_z(f)\geq f^*{\alpha}(z)$ for
all $z\in \D$.  Hence $v\leq v_f^{\alpha}$.  The function
$v_f^{\alpha}$ is harmonic outside the countable set where it takes
the value $-{\infty}$.  Since $v(0)=v_f^{\alpha}(0)$, the maximum
principle implies that $v=v_f^{\alpha}$ in $\D\setminus
(v_f^{\alpha})^{-1}(-{\infty})$, so $v=v_f^{\alpha}$ on $\D$.  We have
proved the following.

\proclaim{2.1. Proposition}  Let ${\alpha}$ be a nonnegative function 
on a complex manifold $X$.  If
$f\in \O(\D,X)$ is an extremal disc for the Lelong
functional $H_L^\alpha$ in the sense that 
$G_{\alpha}(f(0))=H_L^{\alpha}(f)=v_f^{\alpha}(0)$,
and $G_{\alpha}(f(0))>-{\infty}$, then
$$
G_{\alpha}(f(z))= 
\sum\limits_{w\in \D}{\alpha}(f(w)) m_w(f)
\log\bigg|\dfrac{z-w}{1-\bar wz}\bigg|, \qquad z\in \D,
$$
which implies that the function $G_{\alpha}\circ f$ is harmonic outside the 
countable set where it takes the value $-{\infty}$.  
\endproclaim

In \cite{6} we proved that the envelopes of the three functionals are
plurisubharmonic for a large class of manifolds with mild
conditions on ${\varphi}$, $v$, and ${\alpha}$.
We define $\Cal P$ as the class of complex manifolds $X$ for which there
exists a finite sequence of complex manifolds and holomorphic maps 
$$
X_0 @>{h_1}>> X_1 @>{h_2}>>\dots @>{h_m}>> X_m=X,\qquad m \geq 0,
\tag 2.1
$$ 
where $X_0$ is a domain in a Stein manifold and each $h_i$, $i=1,\dots,m$, is
either a covering (unbranched and possibly infinite) or a finite
branched covering (i.e., a proper holomorphic surjection with finite
fibres). 

Assume now that $X\in \Cal P$.  Then $EH_P^\phi$ is plurisubharmonic
for every upper semicontinu\-ous function ${\varphi}$ on $X$
\cite{6, Thms\. 2.2 and 3.4}.

The function $EH_R^v$ is plurisubharmonic for every continuous
plurisubharmonic function $v:X\to\R$ \cite{6, Thm\. 4.4}.

The function $EH_L^\alpha$ is plurisubharmonic for every nonnegative
function ${\alpha}$ on $X$ for which the sequence (2.1) can be chosen
such that ${\alpha}^{-1}[c,\infty)\setminus B$ is dense in
${\alpha}^{-1}[c,\infty)$ in the analytic Zariski topology on $X$ for
every $c>0$, where
$$B=\bigcup_{i=1}^m (h_m\circ\dots\circ h_{i+1})(B_i),$$ 
and $B_i$ denotes the (possibly empty) branch locus of $h_i$ \cite{6,
Thm\. 5.12}.
Observe that this condition holds in particular if ${\alpha}=0$ on
$B$, or if $X$ is a domain in a Stein manifold.  

It turns out that $EH_L^{\alpha}$ is related to the Kobayashi
pseudodistance $\kappa_X$ on $X$. By definition $\kappa_X$ is the
largest pseudodistance on $X$ smaller than or equal to $\delta_X$,
where
$$
\delta_X(x,a)=\inf\set{\varrho_\D(z,w)\, ; \, f(z)=x, f(w)=a \text{ for 
some } f\in \O(\D,X)},
$$
and $\varrho_\D$ denotes the Poincar\'e distance in $\D$,
$$\varrho_\D(z,w)=\tanh^{-1}\bigg\vert\frac{z-w}{1-\bar w z}\bigg\vert.$$
We define $k_X=\log \tanh {\delta}_X$.
By composing the map $f$ in the definition of
$\delta_X$ with an automorphism which sends $0$ to $w$ and 
$z$ to the a point on the positive real axis,  and then replacing it by
$z\mapsto f(z/r)$ with $r>1$ and $r$ close to $1$, we see that
 $$
k_X(x,a)= 
\inf\set{\log t\, ;\, t\in (0,1), f(t)=a, f(0)=x,
\text{ for some } f\in \A_X}. 
 $$
We now define
 $$
k_X^{\alpha}(x)
=\inf\limits_{a\in X} {\alpha}(a) k_X(x,a), \qquad x\in X.
 $$
If $X$ is a domain in a Stein manifold, then
$G_{\alpha}=EH_L^{\alpha}=EH_P^{k_X^{\alpha}}$ \cite{6, Thm\. 5.3}.  

Let $G_a$ be the pluricomplex Green function with a logarithmic pole at
$a$.  Then $G_a\leq k_X(\cdot, a)$.  If $u\in \Cal F_L^{\alpha}$ and
$a\in X$, then ${\nu}_u(a)\geq {\alpha}(a)$, so $u\leq
{\alpha}(a)G_a$, and
$$
G_{\alpha}\leq \inf_{a\in X}{\alpha}(a)G_a\leq k_X^{\alpha}.
$$
The above results may be summarized as follows.

\proclaim{2.2. Theorem} 
Let ${\alpha}$ be a nonnegative function on a domain $X$ in a Stein manifold.  
Then 
$$
G_{\alpha}=EH_L^{\alpha}=EH_P^{k_X^{\alpha}}
=\sup\set{u\in \PSH(X)\, ;\, u\leq \inf_{a\in X}{\alpha}(a)G_a}.
$$
\endproclaim

Our next result shows that $G_\alpha$ has no Monge-Amp\gr ere mass where
it is locally bounded.

\proclaim{2.3. Proposition} Let ${\alpha}$ be a nonnegative function
on a complex manifold $X$.  Then the function $G_{\alpha}$ is maximal
in the open subset of $X$ where it is locally bounded.  
\endproclaim

\demo{Proof}  Let $U$ be a relatively compact domain in the open subset
of $X$ where $G_\alpha$ is locally bounded, and let
$v$ be plurisubharmonic on $X$ such that $v\leq G_\alpha$ on $\d U$. 
Let
$$
w=\cases \max\set{v, G_\alpha} & \text{on }U,\\
G_\alpha  & \text{on } X\setminus U, 
\endcases
$$
Then $w$ is plurisubharmonic on $X$ and $w\in\F_L^\alpha$, so $w\leq
G_\alpha$.  Hence, $v\leq G_\alpha$ on $U$.
\qed\enddemo

Now we turn to the boundary behaviour of the supremum of the
Lelong class.  Let $X$ be a domain in a complex manifold $Y$, and $p$ be
a boundary point of $X$.  Recall that a plurisubharmonic function $v$ on
$X$ is called a {\it strong (plurisubharmonic) barrier} at $p$ if
$\lim\limits_{x\to p}v(x)=0$, and $\sup\limits_{X\setminus V}v<0$ for
every neighbourhood $V$ of $p$ in $Y$. 

A relatively compact domain $X$ in a complex manifold is said to be {\it
B-regular} if every continuous function on the boundary of $X$ extends
to a continuous function on the closure of $X$ which is plurisubharmonic
on $X$.  This notion is due to Sibony \cite{13}.  It is easily seen that
a B-regular domain has a strong barrier at every boundary point, and for
domains in $\C^n$, the converse holds.  For a weaker result on an
arbitrary K\"ahler manifold, see Lemma 3.7.  In $\C^n$, strongly
pseudoconvex domains and smoothly bounded pseudoconvex domains of finite
type are B-regular, and B-regular domains are hyperconvex, but not vice
versa. 

\proclaim{2.4.  Proposition} Let $\alpha$ be a nonnegative function on a
domain $X$ in a complex manifold $Y$, and assume that there exists a
strong plurisubharmonic barrier $v$ at $p\in\d X$.  If some $u$ in
$\F_L^\alpha$ is bounded below on a neighbourhood of $p$, then $G_\alpha$
has limit zero at $p$.  
\endproclaim

Example 3.4 shows that this result may fail if existence of a strong
barrier is replaced by hyperconvexity. 

\demo{Proof} Choose a neighbourhood $V$ of $p$ in $Y$ such that 
$u>\beta\in\R$ in a
neighbourhood of $X\cap \overline V$, and choose $c>0$ such that
$\sup\limits_{X\setminus V} v<\beta/c$. Then $u>cv$ in a neighbourhood of
${\partial}V\cap X$, so the function $w$ defined by
$$
w=\cases \max\set{u, cv} & \text{on }X\cap V,\\
u & \text{on } X\setminus V, 
\endcases
$$
is plurisubharmonic on $X$.  Since $u > \beta$ in a neighbourhood of
$X\cap \overline V$, we have $\nu_u=0=\alpha$ there.
Hence, $w\in \Cal F_L^\alpha$, and we get
$$
\liminf_{x\to p} G_\alpha(x) \geq \liminf_{x\to p} w(x)
\geq\lim_{x\to p} cv(x) = 0.\qed
$$ 
\enddemo

Edigarian \cite{3} has proved that if $X_1$ and $X_2$ are domains in
$\C^{n_1}$ and $ \C^{n_2}$ respectively, and $a=(a_1,a_2)\in
X=X_1\times X_2$, then
$$
G_a(x)=\max\{G_{a_1}(x_1), G_{a_2}(x_2)\},\qquad x=(x_1,x_2)\in X.
$$
This is called the {\it product property} of the pluricomplex Green 
function.  By a modification of Edigarian's proof we get the following
result.

\proclaim{2.5. Theorem} Let $\alpha_1$ and $\alpha_2$ be the
characteristic functions of subsets $A_1$ and $A_2$ of complex
manifolds $X_1$ and $X_2$ respectively, and let $\alpha$ denote the
characteristic function of $A=A_1\times A_2$ on the product manifold
$X=X_1\times X_2$, so
$$
\alpha(x)=\min\set{\alpha_1(x_1),\alpha_2(x_2)}, \qquad
x=(x_1,x_2)\in X.  
$$
Then
$$
EH_L^\alpha(x) = \max\{EH_L^{\alpha_1}(x_1), EH_L^{\alpha_2}(x_2)\}.
$$
If $EH_L^{\alpha_1}$ and $EH_L^{\alpha_2}$ are plurisubharmonic, then
$$
G_\alpha(x) = \max\{G_{\alpha_1}(x_1), G_{\alpha_2}(x_2)\}.
$$
\endproclaim

It is easy to see that the product property fails in general.  Let
$\alpha_1=\chi_{\{0\}}$ and $\alpha_2=2\alpha_1$ on $\D$.  Then $\alpha
=\chi_{\{(0,0)\}}$ on $\D\times\D$, and
$$
EH_L^\alpha(z_1,z_2)=G_{(0,0)}(z_1,z_2)=\max\{\log|z_1|, \log|z_2|\},
$$
but
$$
\max\{EH_L^{\alpha_1}(z_1), EH_L^{\alpha_2}(z_2)\}
=\max\{\log|z_1|, 2\log|z_2|\}.
$$
The latter function is not even a Lelong envelope (although it is
presumably the envelope of a disc functional involving {\it directional}
Lelong numbers).  Both functions have
Lelong number $1$ at $(0,0)$ and $0$ elsewhere.

\proclaim{2.6. Lemma}  Let $X$ be a complex manifold, ${\alpha}$ be a
nonnegative function on $X$, and ${\beta}\in (-{\infty},0)$.  If
$EH_L^{\alpha}(x)<{\beta}$, then there exists $f\in \A_X$ with $f(0)=x$
and
finitely many points $a_1,\dots,a_l$ in $\D\setminus\set{0}$ such that
$$
-{\infty} < \sum_{k=1}^{l} \alpha(f(a_{k}))\,m_{a_{k}}(f)
\,\log|a_{k}|<\beta.   \tag 2.2
$$
\endproclaim

\demo{Proof}  By the definition of the envelope, there exists 
$f\in \A_X$ with $f(0)=x$ such that $H_L^{\alpha}(f)<{\beta}$,
and by the definition of $H_L^{\alpha}$ there are finitely many
points $a_1,\dots,a_l$ in $\D$ such that the right inequality in (2.2)
holds. 

If the sum equals $-{\infty}$ and $f$ is nonconstant, then $a_k=0$ for
some $k$ and ${\alpha}(f(0))>0$.  Then we may assume that $l=1$ and $a_1=0$. 
We choose $a \in \D\setminus\set{0}$ so close to $0$ that $f$ is
holomorphic in a neighbourhood of the image of $\overline\D\to \C$,
$z\mapsto z(z-a)$, and $m_0(f)\log |a|<\beta$.  If we replace $a_{1}=0$
by $a$ and $f$ by $z\mapsto f(z(z-a))$, then (2.2) holds.
   
If $f$ is constant, then $m_a(f)=+{\infty}$ for all $a\in \D$ and the
sum in (2.2) equals $-\infty$.  We choose $a\in \D\setminus\set{0}$ so
close to zero that ${\alpha}(x)\log |a| <{\beta}$.  Let $U$ be a 
neighbourhood of $x$ in $X$ with a biholomorphism $\Psi:U\to\D^n$, 
$x\mapsto 0$.  Let $r>1$ and $\Phi=\text{id}\times\Psi:D_r\times U\to 
D_r\times \D^n$, where $D_r=\set{z\in\C \, ;\, |z|<r}$.  Now we set $l=1$,
$a_1=a$ and replace $f$ by the disc 
$$ 
z\mapsto{\pr}\big({\Phi}^{-1}(z,\epsilon z(z-a),0,\dots,0)\big),$$ 
where
${\pr}:\C\times X\to X$ is the projection, and $\epsilon>0$ is
chosen so small that $\epsilon z(z-a)\in\D$ if $z\in \overline \D$. 
Then $f$ is nonconstant, $f(0)=f(a)=x$, and (2.2) holds.  
\qed \enddemo

\demo{Proof of Theorem 2.5} We have
$m_z(f)=\min\{m_z(f_1), m_z(f_2)\}$  for all $f=(f_1,f_2) \in\A_X$. 
Hence
$$
H_L^{\alpha_j}(f_j)=\sum_{z\in\D} \alpha_j(f_j(z))\,m_z(f_j)\,\log|z|
\leq H_L^\alpha(f),\qquad j=1,2,
$$
and 
$$
EH_L^\alpha(x) \geq \max\{EH_L^{\alpha_1}(x_1), EH_L^{\alpha_2}(x_2)\}.
$$

To establish the reverse inequality, we assume that
$EH_L^{\alpha_j}(x_j)<\beta\in(-\infty,0)$ for $j=1,2$, and
show that $EH_L^\alpha(x)<\beta$.   By Lemma 2.6, there are
$f_j\in\A_{X_j}$
with $f_j(0)=x_j$ and 
$a_{jk}\in\D$, $k=1,\dots,l_j$, $j=1,2$, such that 
$$
-{\infty} <
\sum_{k=1}^{l_j} \alpha_j(f_j(a_{jk}))\,m_{a_{jk}}(f_j)
\,\log|a_{jk}|<\beta, \qquad j=1,2.  \tag 2.3
$$
Choose $f_j$ such that $l_j$ is as small as possible. 
Then $\alpha_j(f_j(a_{jk}))=1$ and $a_{jk}\neq 0$ for all $j$ and $k$.

Assume that $|a_{j1}| \leq |a_{j2}| \leq\dots$.  Set
$$
\mu_{jk}=m_{a_{jk}}(f_j), \quad 
\mu_j=\sum_{k=1}^{l_j}\mu_{jk}, \quad
b_j=\prod_{k=1}^{l_j}a_{jk}^{\mu_{jk}}\neq 0, \quad \text{ and } \quad
c_j=a_{jl_j}.
$$  
Then
$$
|c_j|^{\mu_j}e^\beta\leq|b_j|<e^\beta.  \tag 2.4
$$
The second inequality is equivalent to (2.3).
To prove the first one, suppose $|b_j|<|c_j|^{\mu_j}e^\beta$.
Then
$$
\prod_{k=1}^{m_j}\bigg|\frac{a_{jk}}{c_j}\bigg|^{\mu_{jk}}<e^\beta,
$$
where $m_j<l_j$ is the smallest number with $|a_{jk}|=|c_j|$ for
$k>m_j$.  Hence, (2.3) holds with $f_j$ replaced by $z\mapsto
f_j(c_j z)$, $a_{jk}$ replaced by $a_{jk}/c_j$, and $l_j$ replaced by
$m_j$, which contradicts the fact that $l_j$ is minimal.

We define the Blaschke products
$$
B_j(z)=\prod_{k=1}^{l_j}\bigg(\frac{a_{jk}-z}{1-\bar a_{jk}z}\bigg)
^{\mu_{jk}}.
$$
Then $B_j(0)=b_j$.  We may assume that $|b_1|\geq|b_2|$.  
By precomposing $f_1$ by a suitable embedding of $\D$ into $\D$ with 
$0\mapsto 0$, and thereby changing $a_{1k}$ slightly,
we may assume that $B_1(0)$ is not a critical value of $B_1$.
By Schwarz' Lemma, we still have $|b_1|\geq |b_2|$.

We may assume that $b_1=b_2$.  Indeed, if $|b_1|>|b_2|$, choose $t\in(0,1)$
with $t^{-\mu_2}|b_2|=|b_1|$.  Then $|a_{2k}/t|<1$, since by (2.4),
$$
|a_{2k}|^{\mu_2} \leq |c_2|^{\mu_2} \leq |b_2|e^{-\beta} <
|b_2/b_1| = t^{\mu_2}.
$$
Replacing $f_2$ by $z\mapsto f_2(tz)$ and $a_{2k}$ by $a_{2k}/t$, we
get $|b_1|=|b_2|$.  Finally, replacing $f_2$ by 
$z\mapsto f_2(e^{i\theta}z)$, where $e^{i\theta\mu_2}=b_2/b_1$, and
replacing $a_{2k}$ by $e^{-i\theta}a_{2k}$, we get $b_1=b_2$.

Exactly as in \cite{7}, we obtain (possibly infinite)
Blaschke products $\phi_j$ with $0\mapsto 0$ and
$$
\sigma=B_1\circ\phi_1=B_2\circ\phi_2.
$$
Now $|\sigma(0)|<e^\beta$ and $|\sigma|=1$ almost everywhere on $\T$.
Choose $r\in(0,1)$ so close to $1$ that
$$
\log|\sigma(0)|-\frac1{2\pi}\int_0^{2\pi}\log|\sigma(re^{i\theta})|\,
d\theta<\beta.
$$
By the Riesz representation formula, the left hand side equals
$$
\sum_{i=1}^n m_{z_i}(\sigma(r\cdot))\log|z_i|,
$$
where $z_1,\dots,z_n$ are the zeros of $z\mapsto\sigma(rz)$ in $\D$.

Now define $f\in\A_X$ with $f(0)=(x_1,x_2)$ by
$$
f(z)=(f_1\circ\phi_1(rz),f_2\circ\phi_2(rz)),
\qquad z\in\overline\D.
$$
Since $\sigma(rz_i)=0$, we have $\phi_j(rz_i)=a_{jk_j}$ for some
$k_j$, and
$$
m_{z_i}(\sigma(r\cdot))=\mu_{jk_j} m_{z_i}(\phi_j(r\cdot))
=m_{a_{jk_j}}(f_j) m_{z_i}(\phi_j(r\cdot))
=m_{z_i}(f_j\circ\phi_j(r\cdot)).
$$
Since $\alpha(f(z_i))=\min\limits_j\alpha_j(f_j(a_{jk_j}))=1$, we get
$$
H_L^{\alpha}(f)=\sum_{z\in\D} \alpha(f(z))m_z(f)\log|z|
\leq\sum_{i=1}^n m_{z_i}(f)\log|z_i|
=\sum_{i=1}^n m_{z_i}(\sigma(r\cdot))\log|z_i|<\beta,
$$
so $EH_L^\alpha(x)<\beta$.
\qed\enddemo

\vfill\eject

\specialhead 3.  The Green function with a pole along a complex subspace
\endspecialhead

\noindent Let $X$ be a complex manifold and $A$ be a (closed) complex
subspace of $X$, which is the same thing as a coherent sheaf $\Cal
I=\Cal I_A$ of ideals in the sheaf $\O_X$ of holomorphic functions on
$X$.  Suppose the stalk $\Cal I_p$ of $\Cal I$ at $p\in X$ is generated
by germs $f_1,\dots,f_m$.  The plurisubharmonic functions
$\max\limits_{i=1,\dots,m}\log|f_i|$ and $\log\sum\limits_{i=1}^m|f_i|$
have the same Lelong number $\nu_A(p)$ at $p$.  This number is
independent of the choice of generators.  Namely, say $g_1,\dots,g_k$
also generate $\Cal I_p$.  Then $f_i=\sum h_{ij}g_j$ for some
$h_{ij}\in\O_{X,p}$, so $\sum|f_i|\leq c\sum|g_j|$ on a neighbourhood of
$p$ for some constant $c>0$.  Hence, the Lelong number of
$\log\sum|f_i|$ at $p$ is no smaller than that of $\log\sum|g_j|$. 
Interchanging $\set{f_i}$ and $\set{g_j}$, we see that these Lelong
numbers are the same. 

We might call $\nu_A(p)$ the {\it multiplicity} of $A$ at $p$.  If $A$
is smooth at $p$, then $\nu_A(p)=1$, so if $A$ is a submanifold of $X$,
then $\nu_A$ is the characteristic function $\chi_A$ of $A$.  We set
$$\Cal F_A=\set{u\in \PSH(X)\, ;\, u\leq 0, \nu_u\geq \nu_A}, \qquad
G_A=\sup\Cal F_A, \qquad k_A=k_X^{\nu_A}.$$ 
We call $G_A$ the {\it pluricomplex Green function with a logarithmic 
pole along} $A$. 

\smallskip
If $\Cal I$ is principal, meaning that each stalk $\Cal I_p$ is a
principal ideal in $\O_{X,p}$, and $\Cal I$ is neither $0$ nor $\O_X$,
then $A$ is a hypersurface, i.e., of pure codimension 1.  If $\Cal I$ is
reduced, i.e., $A$ is a subvariety of $X$, then the converse holds.  A
principal coherent ideal sheaf different from the zero sheaf is nothing
but an effective divisor.  

Suppose now that $A$ is an effective divisor.  Then the current $[A]$ of
integration over $A$ is a closed positive (1,1)-current on $X$, locally
defined as $\frac 1 \pi \Cal L(\log|h|)$, where $h$ is a local generator
for $\Cal I$.  If $f\in\O(\D,X)$, then there is a positive Borel measure
$f^*[A]$ on $\D$ defined locally as 
$\frac 1 {2\pi}\Delta\log|h\circ f|$, unless 
$h\circ f=0$, in which case we set $f^*[A]=0$.
There is a {\it generalized Riesz functional} $H_R^A$ associated to $A$,
defined by the formula 
$$H_R^A(f)= \int_\D\log|\cdot|\,f^*[A], \qquad
f\in\O(\D,X).$$ 
If $A$ is principal, say $A$ is the divisor of a holomorphic function
$h$ on $X$, then $H_R^A$ is the Riesz functional $H_R^{\log|h|}$.

\remark{\bf 3.1.  Example} Let $X$ be the unit ball in $\C^n$.  We have
$\kappa_X=\delta_X=c_X$, where $c_X$ is the Carath\ac eodory distance on
$X$, and $G_a=k_X(\cdot, a)=\log\tanh\delta_X(\cdot,a)$, $a\in X$.  This
is in fact true on any bounded convex domain in $\C^n$ by work
of Lempert \cite{9}.  Hence, for every submanifold $A$ of $X$, 
$$G_A\leq k_A=\inf_{a\in A}G_a=\log\tanh c_X(\cdot,A)
=\log\tanh\kappa_X(\cdot, A). $$

Now let $A$ be
the hypersurface defined by the equation $z_1=0$, and write
$z\in \C^n$ as $z=(z_1,z')$ with $z_1\in \C$ and $z'\in \C^{n-1}$. Then
$$
G_A(z)=\log\dfrac{|z_1|}{\sqrt{1-|z'|^2}}, \qquad z\in X.
$$
Namely, if $u$ denotes the function on the right, then $u\in\F_A$, so
$u\leq G_A$.  Also, on each disc $D_c=\set{z\in X\, ;\,z'=c}$ where
$c$ is a constant, $u$ is subharmonic with $\Delta u=2\pi\delta_{(0,c)}$
and $u|\d D_c=0$, and $G_A\leq 0$ is subharmonic with $\Delta
G_A\geq 2\pi\delta_{(0,c)}$.  Hence, by the Riesz representation formula,
$G_A\leq u$ on $D_c$.  

For $a\in X$, the pluricomplex Green function 
on $X$ with a logarithmic pole at
$a$ is $G_a=\log|T_a|$, where $T_a$ is any automorphism
of $X$ with $a\mapsto 0$.  For $a\neq 0$,
one such automorphism is given by the formula
$$
T_a(z)=\frac{a-P_a(z)-s_a Q_a(z)}{1-\langle z,a\rangle},
$$
where $P_a$ is the orthogonal projection onto the linear space spanned by
$a$, $Q_a$ is the orthogonal projection onto its orthogonal complement,
and $s_a=\sqrt{1-|a|^2}$. 
If $a=(0,a')$ and $z'=a'$, then $P_a(z)
=(0,z')=a$ and $Q_a(z)=(z_1,0)$, so
$$
T_a(z)=(\frac{-z_1}{\sqrt{1-|z'|^2}},0),
$$
and
$$
G_A(z)=G_{(0,z')}(z)=\inf_{a\in A} G_a(z).
$$
This equality is in fact very exceptional, as Example 3.5 will
indicate.
\endremark

\proclaim{3.2. Proposition}  Let $A$ be an effective divisor in a complex 
manifold $X$.  If $f$ is a holomorphic function generating the ideal
sheaf of $A$ on an open set $U$, then the plurisubharmonic function
$G_A-\log|f|$ on $U\setminus A$ is locally bounded above on $A\cap U$,
so it extends to a plurisubharmonic function on $U$.  Hence,
$$\Cal L(G_A)\geq\pi[A].$$
\endproclaim

\demo{Proof} Let $p\in A\cap U$.  We will show that $G_A-\log|f|$ is
locally bounded above at $p$.  We may assume that $p$ is a smooth point
of the reduction of $A$, since plurisubharmonic functions always extend
across subvarieties of codimension at least 2.  By applying a local
biholomorphism, we may assume that $p=0\in\C^n$, that the unit polydisc
$P$ centred at $p$ is in $U$, and that the reduction of $A$ is given by
the equation $z_1=0$ in $P$.  If $z=(z_1,z')\in P$, then the analytic
disc $\zeta\mapsto(\zeta,z')$ maps $0$ to $(0,z')$ and $z_1$ to $z$, so
$$k_X(z,(0,z')) = \log\tanh\delta_X(z,(0,z')) \leq
\log\tanh\varrho_\D(0,z_1) = \log|z_1|,$$ and $$ G_A(z) \leq k_A(z)\leq
\nu_A(0,z')k_X(z,(0,z'))\leq\nu_A(0,z')\log|z_1|.$$ At $p$, the germ of
$f$ is the product of $z_1^s$ and a unit, where
$s=\nu_A(p)=\nu_{\log|f|}(p)$ is the order of the zero of $f$ at $p$,
which equals the order of vanishing of $f$ along the reduction of $A$ at
$p$, and this is the same at every point of $A\cap P$.  Hence, for
$z=(z_1,z')\in P$, we have $\nu_A(0,z')=s$, so $$G_A(z)-\log|f(z)|\leq
s\log|z_1|-\log|f(z)|,$$ and this is bounded above near $p$.  
\qed\enddemo

\proclaim{3.3.  Theorem} Let $X$ be a relatively compact domain in a
Stein manifold $Y$, and let $A$ be the intersection with $X$ of a
complex subspace $B$ of $Y$.  Then $G_A$ is locally bounded and maximal
on $X\setminus A$,
$$
\nu_{G_A}=\nu_A,
$$
and
$$\align
G_A&=EH_L^{\nu_A} =EH_P^{k_A} =\sup\set{u\in\PSH(X)\,;\, u\leq k_A} \\
&=\sup\set{u\in\PSH(X)\,;\, u\leq\inf_{a\in A}\nu_A(a)G_a}.
\endalign $$
If $X$ has a strong plurisubharmonic barrier at $p\in\d X\setminus B$,
then
$$\liminf_{x\to p} G_A(x)=0.$$
If $A$ is a divisor, then 
$$G_A=EH_R^A=\sup\set{u\in\PSH(X)\,;\, u\leq 0, \Cal L(u)\geq\pi[A]}.$$ 
\endproclaim

The hypotheses of the theorem are satisfied when $X$ is a smoothly
bounded B-regular domain in $\C^n$ and $A$ is the intersection with $X$
of a complex subspace of a neighbourhood of $\overline X$, because
$\overline X$ has a Stein neighbourhood basis \cite{13}.

\demo{Proof} Since $Y$ is Stein, each stalk of the ideal sheaf $\Cal
I_B$ is generated by global sections of $\Cal I_B$ by Cartan's Theorem
A.  Since $X$ is relatively compact in $Y$, there are finitely many
holomorphic functions $f_1,\dots,f_m\in\Cal I_B(Y)$ which generate all
the stalks $\Cal I_{B,p}$, $p\in X$.  We may assume that $|f_i|<1$ on
$X$.  Let $u=\max\limits_{i=1,\dots,m}\log|f_i|$ on $X$.  Then
$u\in\PSH(X)$, $u\leq 0$, and $\nu_u=\nu_A$.  In particular, $u\in\F_A$,
so $u\leq G_A$.  This shows that $G_A$ is locally bounded in $X\setminus
A$, and hence maximal there by Proposition 2.3, and $\nu_{G_A}=\nu_A$. 

The next four equations follow from Theorem 2.2.  The statement
about boundary limits follows from Proposition 2.4.

Let $f\in\A_X$, and suppose $f(0)\notin A$, so $f^{-1}(A)$ is
finite.  Let $b\in\D$, let $h$ be a local generator for $\Cal I_A$ on a
neighbourhood $U$ of $f(b)$, and let $V\subset\D$ be a neighbourhood of
$b$ such that $f(V)\subset U$.  If $z\in V$ and $h(f(z))=0$, then
$m_z(h\circ f)$ is no smaller than $m_z(f)$ times the order of the zero 
of $h$ at $f(z)$, which is the Lelong number $\nu_A(f(z))$ of $\log|h|$ 
at $f(z)$.  Hence, we have
$$f^*[A]=\frac 1{2\pi}\Delta\log|h\circ f| 
= \sum_{z\in V\cap f^{-1}(A)} m_z(h\circ f)\delta_z \geq \sum_{z\in V} 
\nu_A(f(z))m_z(f)\delta_z$$
on $V$.  Integrating the function $\log|\cdot|$ over any measurable
subset $S$ of $V$ with respect to
these measures gives
$$\int_S\log|\cdot|\,f^*[A]\leq \sum_{z\in
S}\nu_A(f(z))m_z(f)\log|z|.$$
This shows that $H_R^A(f)\leq H_L^{\nu_A}(f)$, so $G_A=EH_L^{\nu_A}\geq 
EH_R^A$ on $X\setminus A$.  On $A$, both envelopes equal $-\infty$.

Let 
$$\F_R^A=\set{u\in\PSH(X)\,;\, u\leq 0, \Cal L(u)\geq\pi[A]}.$$
By Proposition 3.2, $G_A\in\F_R^A$, so $G_A\leq\sup\F_R^A$.  

Now take $u\in\F_R^A$ and
$f\in\A_X$.  If $h$ is a local generator for $\Cal I_A$, then $\Cal
L(\log|h|) =\pi[A]\leq\Cal L(u)$, so $2\pi f^*[A]=\Delta\log|h\circ f|\leq
\Delta(u\circ f)$.  Hence, by the Riesz representation formula, 
$$u(f(0))=\frac 1 {2\pi}\int_\T u\circ f\,d\lambda + \frac 1 {2\pi}
\int_\D \log|\cdot|\,\Delta(u\circ f) \leq 
\int_\D\log|\cdot|\,f^*[A] = H_R^A(f).$$
This shows that $\sup\F_R^A\leq EH_R^A$, so
$$G_A\leq\sup\F_R^A\leq EH_R^A\leq EH_L^{\nu_A}=G_A.\qed$$  
\enddemo

\remark{\bf 3.4.  Example} Let $X=\D\times\D$ be the unit bidisc in
$\C^2$.  Then $X$ is hyperconvex but not B-regular.  Let
$A=\set{z_1=0}$.  Let $\alpha_1=\chi_{\set{0}}$ and $\alpha_2=\chi_\D$
on $\D$.  Then
$$\alpha(z)=\min\set{\alpha_1(z_1),\alpha_2(z_2)}=\chi_A(z)$$ for
$z=(z_1,z_2)\in X$.  By Theorem 2.5, or as in Example 3.1,
$$G_A(z)=G_\alpha(z)=\max\set{G_{\alpha_1}(z_1),G_{\alpha_2}(z_2)}
=\max\set{\log|z_1|, -\infty}=\log|z_1|.$$ Clearly, $G_A$ does not go to
zero at all points of $\partial X\setminus A$.  
\endremark

\remark{\bf 3.5.  Example} Let $X$ be the unit ball in $\C^n$.  It is
well known that for $a\in X$, the Lelong functional
$H_L^{\chi_{\set{a}}}$, whose envelope is the pluricomplex Green
function $G_a$ with a logarithmic pole at $a$, has essentially unique
extremal discs, whose images are complex geodesics in $X$.  More
precisely, for $x\in X$, $x\neq a$, there is $f\in\A_X$ with $f(0)=x$,
unique modulo precomposition by a rotation, such that
$G_a(x)=H_L^{\chi_{\set{a}}}(f)$, namely $f(z)=T(z,0,\dots,0)$, where $T$
is an automorphism of $X$ with $T(0)=x$ and
$T^{-1}(a)\in\D\times\set{0}^{n-1}$. 

Now let $A$ be a submanifold of $X$ and $x\not\in A$.  Then 
$\inf\limits_{a\in A}G_a(x)$ is actually a minimum, and 
$G_A(x)=\inf\limits_{a\in A}G_a(x)$ if and
only if there is $b\in A$ and an extremal disc $f\in\A_X$ with $f(0)=x$
such that
$$\log|f^{-1}(b)|=G_b(x)=G_A(x)\leq H_L^{\chi_A}(f)
=\sum_{f(z)\in A}\log|z|.$$
This implies that
$$f(\D)\cap A=\set{b}.$$
In other words, if $G_A(x)=\inf\limits_{a\in A}G_a(x)$, then the complex
geodesic realizing the hyperbolic distance from $x$ to $A$ intersects
$A$ in only one point.

There is no shortage of counterexamples to this.  When $n=2$, take for
instance the smooth curve
$$A=X\cap\set{(z,w)\in\C^2 \, ;\, z^2+w^2=c}, \qquad 0<|c|<1,$$  
which is connected and intersects each complex geodesic through the 
origin in either zero or two points.

This proves the following.  There is a smooth curve in $\C^2$
whose intersection $A$ with the unit ball $X$ is nonempty and connected 
such that:
\roster 
\item $G_A\neq\inf\limits_{a\in A} G_a$, 
\item the functions
$\inf\limits_{a\in A}G_a$, $k_A$, $\log\tanh c_X(\cdot, A)$, and
$\log\tanh\kappa_X(\cdot,A)$ (which are in fact all equal) are not
plurisubharmonic on $X$, and 
\item $\log\tanh c_X(\cdot,A)$ is not dominated by $G_A$,
even though $\log\tanh c_X(\cdot, a)\leq G_a$ for every $a\in X$. 
\endroster
Furthermore, the domain $Y=X\setminus A$ has the same Carath\ac eodory
distance as $X$ because $A$ is removable for bounded holomorphic
functions.  The strongly convex part of $\d Y$ is infinitely distant
from any point of $Y$, so $c_X(y,A)$ is in fact the Carath\ac eodory
distance from $y\in Y$ to $\d Y$.  Hence, $Y$ is an example
of a bounded pseudoconvex domain in $\C^n$ such that $\log\tanh$ of the
Carath\ac eodory distance to the boundary is not plurisubharmonic. 
\endremark

\smallskip Finally, we shall establish a uniqueness property and a
continuity property of the pluricomplex Green function with a
logarithmic pole along a complex subspace $A$.  Although we expect these
results to hold in general, at present we can only prove them when $A$
is a principal divisor.

Recall that if $X$ is a Stein space with $H^2(X,\Bbb Z)=0$, then the
second Cousin problem can be solved on $X$, so every divisor on $X$ is
principal. 

\proclaim{3.6. Uniqueness Theorem}  Let $X$ be a relatively compact
domain in a complex manifold.  Let $A$ be the divisor of a 
holomorphic function $f$ on $X$.  If $u$ is a negative plurisubharmonic
function on $X$ such that
\roster
\item $u$ is locally bounded and maximal in $X\setminus A$,
\item for every $\epsilon>0$ there is a compact subset $K$ of $X$ such
that $G_A\leq u+\epsilon$ on $X\setminus K$, and
\item every point in $A$ has a
neighbourhood $U$ with a constant $C>0$ such that
$$
\log|f|-C\leq u \leq \log|f|+C \qquad\text{on } U,
$$
\endroster
then $u=G_A$.
\endproclaim

\demo{Proof} By Proposition 3.2, $G_A-\log|f|$ extends to a
plurisubharmonic function $\tilde G_A$ on $X$.  The function $\tilde
u=u-\log|f|$ is plurisubharmonic on $X\setminus A$ and locally bounded
on $X$, so it extends to a locally bounded plurisubharmonic function on
$X$.  Since $u$ is maximal on $X\setminus A$, so is $\tilde u$.  Since
$A$ is pluripolar, $\tilde u$ is maximal on $X$ \cite{5, Prop\.  4.6.4},
so by \therosteritem2, $\tilde G_A\leq\tilde u$ and $G_A\leq u$. 
Finally, by \therosteritem3, $u\in\F_A$, so $u=G_A$.  \qed\enddemo

The following lemma solves the Dirichlet problem for the Monge-Amp\gr ere 
operator on a K\"ahler manifold, without continuity of the solution.
On a Stein manifold, the solution is continuous by Lemma 3.8.

\proclaim{3.7. Lemma} Let $X$ be a relatively compact domain in a K\"ahler
manifold (e.g\. a Stein manifold) with a strong plurisubharmonic barrier
at every boundary point.  Let $\phi:\d X\to\Bbb R$ be a continuous
function.  Then there is a unique maximal plurisubharmonic function $u$ 
on $X$ such that
$$\lim_{x\to p} u(x) = \phi(p) \qquad\text{for every } p\in\d X,
$$
namely $u=\sup\F$, where
$$
\F=\set{ v\in\PSH(X) \,;\, v^*|\d X \leq \phi}.
$$
\endproclaim

Here, $v^*$ denotes the upper semicontinuous function 
$p\mapsto\limsup\limits_{x\to p}v(x)$ on $\overline X$.

The K\"ahler condition provides a link between pluripotential theory 
and real potential theory.  It implies that the Laplacian is the trace
of the Levi form \cite{16, p\. 90}, so plurisubharmonic functions are
subharmonic with respect to the associated Riemannian metric.  Here, 
this has the important consequence that the class $\set{v\in\F\, ; \,
v\geq\min\phi}$ is compact.  We do not know if this is true without 
the K\"ahler condition. 

\demo{Proof} By hypothesis, $X$ is regular with respect to the
Laplacian, so there is a continuous function $h$ on $\overline X$,
harmonic on $X$, with $h|\d X=\phi$.  Then $\F$ coincides with the class
$\set{ v\in\PSH(X) \,;\, v\leq h}$.  Hence, $u=\sup\F\leq h$, so
$u^*\leq h$ by continuity of $h$.  Since $u^*$ is plurisubharmonic,
$u^*\in\F$, so $u^*=u$.  This shows that $u\in\F$. 

Let $U$ be a relatively compact domain in $X$, and let
$v$ be plurisubharmonic on $X$ such that $v\leq u$ on $\d U$. 
Let
$$
w=\cases \max\set{u, v} & \text{on }U,\\
u  & \text{on } X\setminus U, 
\endcases
$$
Then $w$ is plurisubharmonic on $X$ and $w^*|\d X=u^*|\d X\leq\phi$, 
so $w\leq u$.  Hence, $v\leq u$ on $U$.  This shows that $u$ is maximal,
and uniqueness follows.

Now let $p\in\d X$ and $\beta<\phi(p)$.  Let $w$ be a strong barrier at $p$. 
Choose a neighbourhood $V$ of $p$ such that $\phi>\beta$ in $\d X\cap
V$, and choose $c>0$ such that $\beta+c\sup\limits_{X\setminus V}
w<\min\phi$.  Then $v=\beta+cw\in\F$, so $v\leq u$, and
$$\beta = \lim_{x\to p} v(x)\leq\liminf_{x\to p} u(x)\leq\phi(p).$$
Since $\beta$ is arbitrary, this shows that $\lim\limits_{x\to p} u(x)=
\phi(p)$.
\qed\enddemo

The next lemma was proved by Walsh \cite{15} for domains in $\C^n$,
and generalized to Banach spaces by Lelong \cite{8}.  Our proof is based
on Lelong's argument.

In a  metric space $(Y,d)$, we let $B(x,{\varrho})$ 
denote the open ball  with centre $x$ and radius
${\varrho}$, and for any $A\subset Y$ we let 
$A_{\varrho}=\set{x\in A\,;\, d(x,{\partial}A)> {\varrho}}$.

\proclaim{3.8. Lemma}  Let $X$ be a relatively compact domain in 
a manifold $Y$.  Assume that there exists a metric defining the 
topology on $Y$ such that for every 
$v\in\PSH(X)$ and $r>0$, there exists a continuous plurisubharmonic
function $v_r:X_r\to\R$ with  
$$
v(x)\leq v_r(x)\leq \sup\limits_{\overline B(x,r)}v,
\qquad  x\in X_r.
$$   
This holds in particular if $Y$ is Stein.
Let $h:\overline X\to\R$ be continuous, and
$$u=\sup\set{v\in \PSH(X)\, ;\, v\leq h}.$$
If $u^*$ is continuous on ${\partial}X$, then $u$ is continuous on $X$.

\endproclaim

\demo{Proof}  Since $h$ is continuous on the compact set $\overline X$,
the function $u$ is plurisubharmonic and real-valued. 
It suffices to prove that for every  
$a\in X$ and ${\varepsilon}>0$, 
there exists a neighbourhood $U$ of $a$ such that
$u(a)-{\varepsilon} <u$ on $U$.

Since $h$ is uniformly continuous on 
the compact set $\overline X$, and $u^*$ is continuous on the 
compact set ${\partial}X$, there is ${\delta}>0$ such that $a\in X_\delta$, 
$$
\gather
|h(x)-h(y)|< \tfrac 12 {\varepsilon} \quad \text{ for all }
(x,y)\in \overline X\times \overline X \text{ with }
d(x,y)<{\delta},\text{ and }
\tag 3.1 \\
|u^*(x)-u^*(y)|<\tfrac 14 {\varepsilon} \quad \text{for all }
(x,y)\in {\partial} X\times \overline X \text{ with }
d(x,y)< 3{\delta}.
\tag 3.2
\endgather
$$
By assumption, there exists a continuous plurisubharmonic $v:X_\delta\to\R$
such that $u(x)\leq v(x)\leq \sup\limits_{\overline B(x,{\delta})}u$ for 
$x\in X_{\delta}$.  We define $w$ on $X$ by 
$$
w=\cases \max\set{v-\frac 12{\varepsilon},u} &\text{on } X_{\delta},\\
u &\text{on } X\setminus X_{\delta}. \endcases
$$
If $x\in X_{{\delta}}\setminus X_{2\delta}$, then
${\delta}< d(x,{\partial}X)\leq 2{\delta}$.  There exists $b\in
{\partial} X$ such that $d(b,x)\leq 2{\delta}$, so for every
$y\in B(x,{\delta})$ we have $d(b,y)<3{\delta}$, and (3.2) implies that
$$
|u(x)-u(y)|<|u^*(b)-u(x)|+|u^*(b)-u(y)|<\tfrac 12 {\varepsilon},
$$
so
$$
v(x)-\tfrac 12{\varepsilon}\leq \sup_{\overline B(x,{\delta})} u
-\tfrac 12{\varepsilon}<u(x).
$$
Hence, $w=u$ on $X\setminus X_{2{\delta}}$, so $w\in \PSH(X)$.   
By (3.1),
$$
v(x)-\tfrac 12{\varepsilon}
\leq \sup_{\overline B(x,{\delta})}u-\tfrac 12{\varepsilon}
\leq \sup_{B(x,{\delta})}h-\tfrac 12{\varepsilon} <h(x)
,\qquad x\in X_\delta,
$$
so $w\leq h$.  Hence, $v-\tfrac 12{\varepsilon}\leq w\leq u$ 
on $X_{\delta}$.
Now there exists a neighbourhood $U$ of $a$ in $X_\delta$ such that 
$v(a)<v(x)+\tfrac 12{\varepsilon}$ for all $x\in U$, and then
$$
u(a)-{\varepsilon}\leq v(a)-{\varepsilon}
<v(x)-\tfrac 12 {\varepsilon}
\leq u(x), \qquad x\in U. 
$$

Assume now that $Y$ is Stein. 
We may assume that $Y$ is a closed submanifold of $\C^N$.
Let $\C^N$ have the euclidean metric and $Y$ have the induced metric.
By Siu \cite{14, Main Thm\. and Cor\. 1},
there is a Stein neighbourhood $V$ of $Y$ in $\C^N$
and a holomorphic retraction ${\sigma}:V\to Y$.  
Let  $v\in \PSH(X)$ and $r>0$.  
Then  $\tilde v =v\circ {\sigma}$ is plurisubharmonic
on $W={\sigma}^{-1}(X)$.  Since 
$X$ is relatively compact, there is ${\varrho}\in (0,r)$ such
that ${\sigma}(\overline B(x,{\varrho}))\subset B(x,r)$ 
for all $x\in X$.  We choose a nonnegative radially symmetric 
${\chi}\in C_0^{\infty}(\C^N)$ with Lebesgue integral $1$
and support in  $B(0,{\varrho})$.
Then the convolution $\tilde v*{\chi}$ 
defines a smooth  plurisubharmonic function on
$W_{\varrho}$. Now take  $v_r=\tilde v*{\chi}|X_r$.  
\qed
\enddemo

Being plurisubharmonic, $G_A$ is quasi-continuous \cite{5, Thm\. 
3.5.5}, but a quasi-continuous function may be discontinuous everywhere. 

\proclaim{3.9.  Theorem} Let $X$ be a relatively compact domain in a
Stein manifold with a strong
plurisubharmonic barrier at every boundary point.  Let $A$ be the
divisor of a holomorphic function $f$ on $X$ which extends to a
continuous function on $\overline X$.  Then the set of points in $X$ at
which $G_A$ is discontinuous is pluripolar.  
\endproclaim

We are unable to prove that $G_A$ is continuous, nor do we have 
counterexamples.  Continuity of the single-pole Green function on a bounded 
hyperconvex domain in a Stein manifold was proved by Demailly \cite{1}. 
Continuity of the Green function with finitely many weighted poles on a 
bounded hyperconvex domain in a Banach space was proved by Lelong \cite{8}. 

\demo{Proof} By Proposition 3.2, $G_A-\log|f|$ extends to a
plurisubharmonic function $\tilde G_A$ on $X$.  By Lemmas 3.7 and 3.8, 
there are continuous
maximal plurisubharmonic functions $v_j$, $j\in\N$, on $X$ such that 
$$\lim_{x\to p}v_j(x)=\min\set{-\log|f(p)|, j} \qquad\text{for
every } p\in\d X.$$
Then $v_j+\log|f|\in\F_A$, so $v_j\leq\tilde G_A$, and the increasing 
sequence $(v_j)$ is locally bounded above.  Hence, $(v_j)$ converges
pointwise to a lower semicontinuous function $v:X\to\R$, whose
upper semicontinuous regularization $v^*$ is plurisubharmonic on $X$,
and $v^*\leq\tilde G_A$.  The subset $N$ of $X$ where $v<v^*$ is
pluripolar.  Furthermore, since the $v_j$ are maximal, so is $v^*$
\cite{5, Thm\. 3.6.1}.

For $p\in\d X$, we have
$$\liminf_{x\to p} v(x) \geq \liminf_{x\to p} v_j(x)
=\min\set{-\log|f(p)|, j} $$
for each $j$, so letting $j\to\infty$ we get
$$\limsup_{x\to p}\tilde G_A(x)\leq -\log|f(p)|\leq\liminf_{x\to p}v(x)
\leq \liminf_{x\to p}v^*(x).$$
Since $v^*$ is maximal, this shows that $\tilde G_A\leq v^*$, so $\tilde
G_A = v^*$.  

For $p\in X$, 
$$v(p)=\liminf_{x\to p} v(x) \leq \liminf_{x\to p} \tilde G_A(x)
\leq \limsup_{x\to p}\tilde G_A(x) = \tilde G_A(p) = v^*(p),$$
so $\tilde G_A$ is continuous at $p$ if $p\notin N$.  Hence, $G_A$ is
continuous at all points outside the pluripolar set $N\setminus A$.
(Note that this is stronger than saying that the restriction of $G_A$ to
the complement of a pluripolar set is continuous.)
\qed\enddemo

\enddocument